\documentclass[11pt,a4paper]{article}
\usepackage[utf8]{inputenc}
\usepackage[english]{babel}
\usepackage{graphicx}
\usepackage{float}
\usepackage{pdfpages}
\usepackage{amsthm}
\usepackage{url}
\usepackage{xcolor}
\usepackage{subcaption}
\usepackage[font=small,labelfont=rm]{caption}
\newtheorem{theorem}{Theorem}[section]
\newtheorem{lemma}{Lemma}[section]

\newtheorem{assumption}{Assumption}[section]

\theoremstyle{definition}
\newtheorem{exmp}{Example}[section]

\theoremstyle{definition}

\usepackage{authblk}
\usepackage{amsmath}
\usepackage{amsfonts}
\usepackage{amssymb}
\usepackage{algorithmic}
\usepackage[linesnumbered,ruled,vlined]{algorithm2e}
\numberwithin{equation}{section}
\usepackage{mathtools}
\usepackage{scalefnt}
\usepackage{empheq}
\usepackage[a4paper,left=3.5cm,right=3cm,top=3cm,bottom=3cm]{geometry}
\usepackage{makeidx}
\usepackage{bm}
\makeindex

\usepackage{indentfirst}

\date{\normalsize October 25, 2025}

\title{An initial-boundary corrected splitting method for diffusion-reaction problems}

\author[1]{Thi Tam Dang\thanks{tam.dang@helsinki.fi}}
\author[2]{Lukas Einkemmer\thanks{lukas.einkemmer@uibk.ac.at}}
\author[2,3]{Alexander Ostermann\thanks{alexander.ostermann@uibk.ac.at}}
\affil[1]{\normalsize Department of Mathematics and Statistics, University of Helsinki, Finland}
\affil[2]{\normalsize Department of Mathematics, University of Innsbruck, Austria}
\affil[3]{\normalsize Digital Science Center, University of Innsbruck, Austria}

\begin{document}
\maketitle	

\begin{abstract}
Strang splitting is a widely used second-order method for solving diffusion-reaction problems. However, its convergence order is often reduced to order $1$ for Dirichlet boundary conditions and to order $1.5$ for Neumann and Robin boundary conditions, leading to lower accuracy and reduced efficiency. In this paper, we consider a new splitting approach, called \emph{an initial-boundary corrected splitting}, which avoids order reduction while improving computational efficiency for a wider range of applications. In contrast to the corrections proposed in the literature, it does not require the computation of correction terms that depend on the boundary conditions and boundary data. Through rigorous analytical convergence analysis and numerical experiments, we demonstrate the improved accuracy and performance of the proposed method.
\end{abstract}

\section{Introduction}	
\label{sec:intro}
Our goal is to develop an effective numerical approach for solving diffusion-reaction problems. This topic has received considerable attention, especially in the context of time integration schemes. Splitting methods, in particular, have received significant interest (see, e.g., references \cite{f632bc39-61b5-3c7b-90cc-3a08f6b7fde1,GERISCH2002159,hairer2013geometric,HUNDSDORFER1995191,hundsdorfer2013numerical,mclachlan_quispel_2002}).  This is because they allow us to treat the diffusion and reaction problems separately. The (often) linear diffusion problem can be efficiently solved using fast Poisson solvers, such as multigrid or potential methods (see \cite{MCKENNEY1995348}). The (often nonlinear) reaction problem can be solved pointwise. Another favorable property of splitting methods is that they preserve positivity, regardless of the time step size, if the sub-solvers have this property (see \cite{HANSEN20121428}).

Strang splitting~\cite{doi:10.1137/0705041} achieves second-order convergence for diffusion-reaction problems with simple boundary conditions (e.g., periodic boundary conditions). However, in general, second-order convergence requires compatibility between the diffusion boundary conditions and the reaction flow. This is often not the case, in particular, for inhomogeneous Dirichlet or Neumann boundary conditions, resulting in order reduction (see, e.g., the experiments in \cite{doi:10.1137/19M1257081,EINKEMMER201876,doi:10.1137/140994204,doi:10.1137/16M1056250}).

The issue of order reduction has attracted considerable interest from researchers, resulting in several methods to overcome it. In a series of papers Einkemmer and Ostermann \cite{EINKEMMER201876,doi:10.1137/140994204,doi:10.1137/16M1056250,Einkemmer2016ACO} introduced a class of modified splitting methods for diffusion-reaction equations. Their approach modifies both partial flows by adding and subtracting a correction term, respectively, without changing the boundary conditions. Bertoli and Vilmart~\cite{doi:10.1137/19M1257081} proposed a different method in which the correction term is constructed from the output of the flow, rather than directly from the nonlinearity $f$. Finally, Alonso-Mallo, Cano, and Reguera \cite{10.1016/j.cam.2019.02.023} avoid the order reduction by adapting the boundary values in an appropriate way.

Existing modified splitting approaches usually require the precomputation of a correction term, as discussed above. This complicates the implementation and increases the computational complexity. In this paper, we propose an initial-boundary corrected (IBC) splitting scheme for diffusion-reaction problems that overcomes order reduction and simplifies implementation. Subtracting a smooth correction term $z(t)$, constructed from the numerical solution, enforces zero initial data and homogeneous boundary conditions, thereby avoiding the need for precomputed corrections for different boundary conditions. This makes the method more efficient and broadly applicable than existing approaches.

In addition to introducing our novel IBC splitting, we provide a rigorous theoretical analysis of the proposed Strang splitting method. Specifically, Theorem~\ref{Th: C1} establishes second-order convergence of the method within the framework of analytic semigroups under suitable smoothness assumptions. To demonstrate the practical impact of our approach, we present numerical experiments confirming the convergence and illustrating the method's effectiveness compared to existing approaches.

The paper is organized as follows: In Section~\ref{sec:problem}, we introduce our model problem, which is a diffusion-reaction problem with oblique boundary conditions. Section~\ref{sec:splitting} presents the construction of our new splitting approach. Section~\ref{sec:convergence} provides a thorough convergence analysis within a framework of analytic semigroups. We prove that the IBC Strang splitting scheme has a global error of order two. Finally, we conclude our paper by presenting a series of numerical experiments in both one and two dimensions. These experiments serve to illustrate our convergence results and to demonstrate the improved accuracy of our approach.

\section{Problem setting}
\label{sec:problem}	
This section introduces the model problem for which we conduct a thorough convergence analysis of the initial-boundary corrected Strang splitting.  Let $\Omega \subset \mathbb{R}^{d}$ be a bounded open subset with a smooth boundary $\partial \Omega$. We consider the following diffusion-reaction problem in $\Omega$ with oblique time-invariant  boundary conditions
\begin{equation}\label{C2.1}
	\begin{aligned}
	\partial_{t} u(t) &= Du(t) + f(t, u(t)), \\
	Bu|_{\partial \Omega} &= b,  \\
	u(0)&= u_0.	
	\end{aligned}
\end{equation}	
The second-order linear elliptic differential operator $D$ is given by
\begin{equation*}\label{C2.2}
	Du = \sum_{i,j=1}^{d} \partial_{i}(a_{ij}(x) \partial_{j} u)  + \sum_{i=1}^{d} c_{i}(x)\partial_{i}u+ d(x) u,
\end{equation*}
where the matrix-valued function $(a_{ij}(x)) \in \mathbb{R}^{d \times d}$ is assumed to be symmetric, and the coefficients $a_{ij}, c_{i}, d$ are sufficiently smooth. When $a_{ij}(x)= \delta_{ij}$, $b_{i}=0$, $c = d = 0$, the operator $D$ is the Laplacian.

The operator $D$ is called strongly elliptic if, for any $ x \in \Omega$, $\sigma \in \mathbb{R}^{d}$, there exists a positive constant $\eta$ such that
\begin{align*}
	\sum_{i,j=1}^{d} a_{ij}(x) \sigma_{i}\sigma_{j} \ge \eta \sum_{i=1}^{d}\sigma_i^2.
\end{align*}
For $\beta_{i}(x)$ and $\alpha(x)$ be sufficiently smooth, we consider the boundary operator $B$ given by
\begin{equation*}
	Bu = \sum_{i=1}^{d} \beta_i(x)\partial_{i}u + \alpha(x)u.
\end{equation*}
We assume that $B$ satisfies the condition of uniform non-tangency as described in $\cite{doi:10.1137/16M1056250}$. The choice of the coefficients $\beta_{i}$ and $\alpha$ determines the type of boundary conditions applied. For example, setting  $\alpha = 1$ and $ \beta_{i}(x)= \sum_{j=1}^{d} a_{ij}(x) n_{j}(x)$ for $1 \le i \le d$, where $n(x)$ represents the outward normal to $\Omega$ at any point $x \in \partial \Omega$, yields a Robin boundary condition.

Furthermore, the nonlinearity $f$, the boundary data $b$, and the initial condition $u_0$ are assumed to be sufficiently smooth. This implies that $f$ is locally Lipschitz continuous in a neighborhood of the exact solution, which is used in the proof of Theorem \ref{Th: C1}.

\section{Description of the proposed splitting method}
\label{sec:splitting}
In this section, we present the construction of the proposed IBC splitting method for the numerical solution of \eqref{C2.1}. The main idea is to transform the problem \eqref{C2.1} so that the initial data become zero and the boundary data are homogeneous. This is achieved by constructing a piecewise smooth function $z(t)$ that satisfies the boundary conditions while preserving the initial data.

We start with describing how to perform a splitting step in the time interval \([t_{n}, t_{n+1}]\) with step size \(\tau\). Let \( u_n \) be the numerical approximation of \( u(t) \) at time \( t = t_n \) and $z_n(t)$ the desired correction. We simply set \( z_n(t) = u_n \), since the numerical solution \(u_{n}\) satisfies the boundary and the initial condition for the step, i.e.~$z_n(t_n)=u_n$. The transformed variable $\hat{u}_n(t)= u(t)-z_n(t) = u(t)-u_n$ is a solution of the following problem on the time interval \( t \in [t_{n}, t_{n+1}]\):
\begin{subequations}\label{3.1}
    \begin{align}
    \partial_{t} \hat{u}_{n} &= D\hat{u}_{n} + h(t,\hat{u}_{n}) + g_{n}(t), \label{3.1a} \\
    B\hat{u}_{n}|_{\partial \Omega}&= 0, \label{3.1b}
    \end{align}
\end{subequations}
where $h$ is the modified nonlinearity that has the form
\begin{align}\label{C3.2}
    h(t,\hat{u}_{n})= f(t, \hat{u}_{n}+ u_{n}) -f(t, u_{n}),
\end{align}
and $g_n$ is the following source term
\begin{align}\label{3.3}
    g_{n}(t) = D u_n + f(t, u_n).
\end{align}
We now split \eqref{3.1} into two subproblems:
\begin{equation}\label{3.4}
    \begin{aligned}
    \partial_{t} \hat{v}_{n}&=  D \hat{v}_{n} + g_{n}(t),\\
    B\hat{v}_{n}|_{\partial \Omega}&= 0
    \end{aligned}
\end{equation}
and
\begin{equation}\label{3.5}
    \partial_{t} \hat{w}_{n}= h(t, \hat{w}_{n}),
\end{equation}
and apply standard Strang splitting to this set of problems with the initial value $\hat v_n(t_n)=0$. This approach will be called the initial-boundary corrected splitting scheme for time-invariant boundary conditions. The associated numerical scheme is given in Algorithm \ref{alg 1}.

\medskip
\begin{algorithm}[H]
\SetAlgoLined
Solve \eqref{3.4} with $g_n(t)$ given by \eqref{3.3} using the initial data \(\hat{v}_{n}(t_n)=0 \) to obtain \( \hat{v}_{n}(t_n + \frac{\tau}2)\)\;
Solve \eqref{3.5} using the initial data \( \hat{w}_{n}(t_n) = \hat{v}_{n}(t_n+\frac{\tau}2) \) to obtain \( \hat{w}_{n}(t_n+\tau)\)\;
Solve \eqref{3.4} with $g_n(t)$ given by \eqref{3.3} with initial data \( \hat{v}_{n}(t_n+\frac{\tau}2) = \hat{w}_{n}(t_n+ \tau)\) to obtain \(\hat{u}_{n+1} = \hat{v}_{n}(t_n+\tau)\)\;
Set $u_{n+1}= \hat{u}_{n+1} + u_n$\;
\caption{IBC Strang splitting for \eqref{C2.1}; time-invariant boundary data $b$.}
\label{alg 1}
\end{algorithm}

\medskip

We note that the source term \(g_{n}(t)\) given by \eqref{3.3} can be easily obtained from $u_n$ by applying operations that are, in any case, necessary for computing the right-hand side of the equation. We also note that most fast Poisson solvers can easily tolerate the additional source term in \eqref{3.4}.

\section{Convergence analysis}
\label{sec:convergence}
The aim of this section is to provide a comprehensive analysis of the convergence properties of the initial-boundary corrected Strang splitting method applied to \eqref{C2.1}. We consider the abstract evolution equation
\begin{equation}\label{C4.4}
    \partial_{t}\hat{u}(t) + A \hat{u}(t) = h(t,\hat{u}(t))+ g(t),
\end{equation}
where $h$ is defined by \eqref{C3.2} and $-A$ is the infinitesimal generator of an analytic semigroup on the Banach space $X$ with norm $\|\cdot\|$. The action of $A$ is defined by $Aw = -Dw$ for all functions $w\in D(A)$, i.e.~in the domain of $A$. This subspace of $X$ includes the homogeneous boundary conditions associated with the problem. For example, for $X=L^2(\Omega)$ and $D$ being a second-order strongly elliptic operator, we have $D(A) = H^2(\Omega)\cap H^1_0(\Omega)$ for Dirichlet boundary condition.

On bounded time intervals $t\in [0,T]$, analytic semigroups are bounded
\begin{align}\label{C4.1}
	\Vert e^{-tA}\Vert \le C,
\end{align}
and possess the parabolic smoothing property
\begin{align}\label{C4.2}
	\| A^{\gamma} e^{-tA} \| \le C t^{- \gamma} \quad \text{for all}\quad \gamma > 0.
\end{align}
For details, we refer to the textbooks \cite{henry1981geometric} and \cite{pazy2012semigroups}.

Applying the variation of constants formula to \eqref{C4.4} allows us to write the exact solution of \eqref{C2.1} as follows
\begin{equation}\label{C4.6}
\begin{aligned}
u(t_{n+1})&= z_n(t_{n+1})+ e^{-\tau A} \hat{u}(t_{n})+ \int_0^{\tau} e^{-(\tau-s)A} g(t_{n}+s) \,ds \\
& \qquad +\int_0^{\tau}e^{-(\tau-s)A}  h(t_{n}+s, \hat{u}(t_{n}+s)) \,ds.
\end{aligned}
\end{equation}
We now split \eqref{C4.4} into the diffusion equation
\begin{equation}\label{4.7}
\partial_{t} \hat{v}(t) + A \hat{v}(t)= g(t),	
\end{equation}
and the nonlinear reaction equation
\begin{equation}\label{4.8}
\partial_{t} \hat{w}(t) = h(t, \hat{w}(t)).
\end{equation}
The general idea of the convergence analysis follows the approach outlined in \cite{doi:10.1137/16M1056250}. The main differences are due to the imposition of zero initial data on the numerical solution, which allows for the elimination of a number of terms and consequently a simple and concise proof.

\subsection{Error analysis}
\label{subsec:error-anal}
First, we study the local error of the initial-boundary corrected Strang splitting method when applied to \eqref{C2.1}. To do this, we consider one step of the numerical solution starting at time $t_{n}$ on the exact solution. Applying the correction to the initial data gives the transformed initial data $\hat{v}(t_{n})= \hat{u}(t_{n})=u(t_{n})-u_{n}$. The process now starts with the integration of \eqref{4.7} with a step of size $\tau/2$, which gives
\begin{equation}\label{C4.7}
    V = \hat{v}\left(t_{n}+ \tfrac{\tau}{2}\right)= e^{-\frac{\tau}{2}A} \bigl(u(t_{n})-u_{n}\bigr)+\int_0^{\frac{\tau}{2}} e^{-\left(\frac{\tau}{2}-s\right)A } g(t_{n}+s)\,ds.
\end{equation}
The subsequent solution of \eqref{4.8} can be expressed in terms of a Taylor series expansion as follows
\begin{equation}\label{C4.8}
    \hat{w}(t_{n}+\tau) = V + \tau h(t_{n}, V) + \frac{\tau^2}{2} \Big( \partial_1 h(t_{n}, V)+ \partial_2 h(t_{n}, V) h(t_{n}, V) \Big) + \mathcal{O}(\tau^3).
\end{equation}
Here, $\partial_1 h$ and $\partial_2 h$ denote the partial derivatives of $h$ with respect to its first and second arguments, respectively. We write $\mathcal{O}(\tau^3)$ for the remainder term, which is bounded and of order~3 in $\tau$. This notation for the remainder term is used consistently throughout this subsection. We finally integrate \eqref{4.7} with a step of size $\tau/2$, using the initial value $\hat{v}\left( t_{n}+ \frac{\tau}{2}\right) = \hat{w}(t_{n}+\tau) $ to get
\begin{equation}\label{C4.9}
    \hat u_{n+1} = e^{-\frac{\tau}{2}A} \hat{w}(t_{n}+\tau)+\int_0^{\frac{\tau}{2}} e^{-\left(\frac{\tau}{2}-s \right)A } g\bigl(t_{n}+ \tfrac{\tau}{2}+s \bigr)\, ds.
\end{equation}
The numerical solution of \eqref{C2.1} at time $t_{n+1}$ with initial value $u(t_n)$ at time $t_n$ is then given by
\begin{equation}\label{C4.11}
    \begin{aligned}
    \mathcal{S}_{\tau}u(t_{n})& = z_n(t_{n+1})+ \hat u_{n+1}\\
    & = z_n(t_{n+1})+ e^{-\tau A} \bigl(u(t_{n})-u_{n}\bigr)+ \tau e^{-\frac{\tau}{2}A} h(t_{n}, V) +\int_0^{\tau} e^{-(\tau-s)A} g(t_{n}+s)\,ds \\
    & \qquad+ \frac{\tau^2}{2} e^{-\frac{\tau}{2}A} \Big(\partial_1 h(t_{n}, V)+ \partial_2 h(t_{n}, V) h(t_{n}, V)\Big) + \mathcal{O}(\tau^3),
    \end{aligned}
\end{equation}
where we have denoted the IBC Strang splitting operator by $\mathcal{S}_{\tau}$.

Let the local error be denoted by $\delta_{n+1}= \mathcal{S}_{\tau} u(t_{n})-u(t_{n+1})$. Subtracting the exact solution \eqref{C4.6} from the numerical solution obtained in \eqref{C4.11} we obtain
\begin{equation}\label{C4.12}
    \begin{aligned}
    \delta_{n+1}&= \tau e^{-\frac{\tau}{2}A} h(t_{n}, V) + \frac{\tau^2}{2} e^{-\frac{\tau}{2}A} \Big(\partial_1 h(t_{n}, V)+ \partial_2 h(t_{n}, V) h(t_{n}, V)\Big) \\
    & \qquad - \int_0^{\tau} e^{-(\tau -s)A} h\left(t_{n}+s, \hat{u}(t_{n}+s)\right) ds + \mathcal{O}(\tau^3).
    \end{aligned}
\end{equation}
For our analysis, we use the following assumption on the data of \eqref{C2.1}.

\begin{assumption}\label{as1}
Let $\Omega$ be a domain with smooth boundary. Let the linear differential operator $D$ be strongly elliptic with sufficiently smooth coefficients subject to oblique boundary conditions with smooth data. Further let the nonlinearity $f$ and the solution $u$ also be sufficiently smooth.
\end{assumption}

Lemma \ref{Le2} is a refined estimate of the local error~\eqref{C4.12} that we use in the proof of Theorem $4.1$ to show convergence of the global error.

\begin{lemma}\label{Le2}
Under the Assumption \ref{as1}, the initial-boundary corrected Strang splitting applied to \eqref{C2.1} satisfies the local error bounds
\begin{equation}\label{C4.23}
    \delta_{n+1} = \mathcal{O}(\tau^2), \qquad \delta_{n+1} = A \hat{\delta}_{n+1}+ \mathcal{O}(\tau^3), \qquad \hat{\delta}_{n+1}= \mathcal{O}(\tau^3).
\end{equation}
\end{lemma}

\begin{proof}
We set
\begin{equation}\label{C4.14}
    \phi_{n}(s)= e^{-(\tau-s)A} \hat{\phi}_{n}(s), \quad \hat{\phi}_{n}(s)=h( t_{n}+s, \hat{u}(t_{n}+s)).
\end{equation}	
Using the midpoint rule, one gets the following identities
\begin{equation}\label{C4.24}
    \begin{aligned}
    \int_0^{\tau} \phi_{n}(s)\, ds &= \tau \phi_{n}\Bigl( \frac{\tau}{2} \Bigr) + \int_0^{\tau}\int_{\frac{\tau}{2}}^{s} \phi_{n}^{\prime}(\xi)\, d\xi\, ds\\
    &= \tau \phi_{n}\Bigl( \frac{\tau}{2}\Bigr) + \frac{1}{2}\int_0^{\tau} M(s,\tau) \phi_{n}^{\prime \prime}\left(s\right) ds,
    \end{aligned}
\end{equation}
with the kernel $M(s,\tau)= s^2$ if $s< \tau/2$ and $M(s,\tau) = (\tau-s)^2$ if $s>\tau/2$.	Using \eqref{C4.24}, we can rewrite the local error \eqref{C4.12} as follows
\begin{equation}\label{C4.25}
    \delta_{n+1} = \delta_{n+1}^{\left[1 \right] } + \delta_{n+1}^{\left[2 \right] },
\end{equation}
where
\begin{equation}\label{C4.26}
    \begin{aligned}
    \delta_{n+1}^{\left[1 \right] } &= \tau e^{-\frac{\tau}{2}A} h(t_{n}, V) -\tau e^{-\frac{\tau}{2}A} h\left( t_{n}+ \tfrac{\tau}{2}, \hat{u}\left( t_{n}+\tfrac{\tau}{2}\right) \right) \\
    & \qquad + \frac{\tau^2}{2} e^{-\frac{\tau}{2}A} \Big(\partial_1 h(t_{n}, V)+ \partial_2 h(t_{n}, V) h(t_{n}, V)\Big)
    \end{aligned}
\end{equation}
and
\begin{equation}\label{eq:d2}
    \delta_{n+1}^{\left[2 \right] }= -\frac{1}{2} \int_0^{\tau} M(s,\tau) \phi_{n}^{\prime \prime}\left(s\right)ds+ \mathcal{O}(\tau^3).
\end{equation}
Let us first bound $ \| \delta_{n+1}^{\left[1 \right] } \|$. Using \eqref{C4.6} and \eqref{C4.7}, we get
\begin{equation}\label{C4.28}
    \hat{u}\left( t_{n}+ \tfrac{\tau}{2} \right) - V= \int_0^{\frac{\tau}{2}} e^{-\left(\frac{\tau}{2}-s\right) A} h(t_{n}+s, \hat{u}\left(t_{n}+s)\right) ds,
\end{equation}
where the integrand
$
\tilde{\phi}_{n}(s) = e^{-\left(\frac{\tau}{2}-s\right) A} h(t_{n}+s, \hat{u}\left(t_{n}+s)\right)
$
satisfies the identity
\begin{equation}\label{C4.expand}
    \int_0^{\frac{\tau}2} \tilde \phi_{n}(s)\, ds = \frac{\tau}2 h\bigl(t_{n}+\tfrac{\tau}2, \hat{u}\bigl(t_{n}+\tfrac{\tau}2\bigr)\bigr) + \int_0^{\frac{\tau}2}\int_{\frac{\tau}{2}}^{s} \tilde\phi_{n}^{\prime}(\xi)\, d\xi\, ds.
\end{equation}
Using \eqref{C4.28} and \eqref{C4.expand} thus shows that
\begin{equation}\label{C4.29}
    \begin{aligned}
    h\bigl( t_{n}+ \tfrac{\tau}{2}, \hat{u}\bigr( t_{n}&+\tfrac{\tau}{2}\bigr) \bigr)- h(t_{n}, V)\\
    &=  \frac{\tau}{2} \partial_1h(t_{n}, V)+ \frac{\tau}{2} \partial_2 h(t_{n}, V) h(t_{n}, \hat{u}(t_{n}))\\
    & \qquad + \partial_2h(t_{n}, V) \int_0^{\frac{\tau}{2}} \int_{\frac{\tau}{4}}^{s} \tilde{\phi}_{n}^{\prime} (\xi)\,d\xi\, ds + \mathcal{O}(\tau^2).
    \end{aligned}
\end{equation}
Inserting \eqref{C4.29} into \eqref{C4.26} we get
\begin{equation}\label{C4.30}
    \begin{aligned}
    \delta_{n+1}^{\left[1 \right] }& = -\tau e^{-\frac{\tau}{2}A}\partial_2 h(t_{n}, V) \int_0^{\frac{\tau}{2}} \int_{\frac{\tau}{4}}^{s} \tilde{\phi}_{n}^{\prime}(\xi)\, d\xi \,ds + \mathcal{O}(\tau^3)\\
    & = - \tau e^{-\frac{\tau}{2}A}\partial_2 h(t_{n}, V) \int_0^{\frac{\tau}{2}} \int_{\frac{\tau}{4}}^{s} e^{-(\frac{\tau}{2} - \xi)A} \big( A \hat{\phi}_{n}(\xi) + \hat{\phi}_{n}^{\prime}(\xi) \big) d\xi \,ds.
    \end{aligned}
\end{equation}
Since the compatibility condition $h(t,0)=0$ is satisfied we can bound a single application of $A$. Using \eqref{C4.1} and the Assumption \ref{as1} we get
\begin{align}\label{4.21}
    \| e^{-(\frac{\tau}{2} - \xi)A} \big( A \hat{\phi}_{n}(\xi) + \hat{\phi}_{n}^{\prime}(\xi) \big)\| \le C.
\end{align}
This implies that
\begin{equation}\label{C4.31}
    \delta_{n+1}^{\left[1 \right] } = \mathcal{O}(\tau^3).
\end{equation}
It remains to bound $\delta_{n+1}^{[2]}$. Using \eqref{C4.24} and the previous bound for ${\phi}_{n}^{\prime}$ shows the preliminary estimate
$$
\delta_{n+1}^{[2] } = \mathcal{O}(\tau^2).
$$
A refined estimate is obtained by considering \eqref{eq:d2} and bounding $\phi_{n}^{\prime \prime}(s)$ by
\begin{align}\label{C4.32}
    \phi^{\prime \prime}_{n}(s)= A \Big(e^{-(\tau-s)A}A \hat{\phi}_{n}(s)+2e^{-(\tau-s)A} \hat{\phi}_{n}^{\prime}(s)\Big)+ e^{-(\tau-s)A} \hat{\phi}_{n}^{\prime \prime}(s).
\end{align}
Using \eqref{C4.1}, \eqref{4.21}, and the Assumption \ref{as1}, we obtain
\begin{equation}
    \delta_{n+1}^{[2]} = A \hat{\delta}_{n+1} + \mathcal{O}(\tau^3).
\end{equation}
Combining this with \eqref{C4.31} completes the proof of Lemma \ref{Le2}.
\end{proof}

We are now in a position to state the main convergence result in Theorem \ref{Th: C1}.

\begin{theorem}\label{Th: C1}
Suppose that Assumption \ref{as1} holds. Then, there exists a constant \( \tau_0 > 0 \) such that for all step sizes \( 0 < \tau \le \tau_0 \) and for \( t_{n} = n\tau \), the initial-boundary corrected Strang splitting method satisfies the bound
\begin{align}\label{4.25}
	\| u_{n} - u(t_{n}) \| \le C \tau^2 (1 + \left|\log \tau\right|), \quad 0 \le n\tau \le T,
\end{align}
where the constant \( C \) depends on \( T \) but is independent of \( \tau \) and \( n \).
\end{theorem}

\begin{proof}
Let $e_{n}=u_{n}-u(t_{n})$ be the global error. We express $e_{n}$ as follows
\begin{align}\label{C4.35}
    e_{n+1}=\mathcal{S}_{\tau}u_{n}-\mathcal{S}_{\tau} u(t_{n})+\delta_{n+1},
\end{align}	
where $\delta_{n+1}= \mathcal{S}_{\tau} u(t_{n})- u(t_{n+1})$ denotes the local error. Similar to \eqref{C4.11}, we first solve  \eqref{3.4} with the initial value \(\hat{v}_{n} =0 \) to get the following expression
$$
V_{n} = \hat{v}_{n}\bigl(t_n+\tfrac{\tau}2\bigr) =  \int_0^{\frac{\tau}{2}} e^{ -\left( \frac{\tau}{2}-s \right)A } g(t_{n}+s)\, ds.
$$
We note for later use that $\|V-V_n\| \le C \|e_n\|$. Next, we perform the steps in Algorithm~\ref{alg 1} to obtain the numerical solution
\begin{equation}\label{C4.39}
    \begin{aligned}
    \mathcal{S}_{\tau}u_{n}& = z_n(t_{n+1})+ \tau e^{-\frac{\tau}{2}A} h(t_{n}, V_{n}) + \int_0^{\tau} e^{-(\tau-s)A} g(t_{n}+s)\, ds \\
    & \qquad + \frac{\tau^2}{2} e^{-\frac{\tau}{2} A} \Big( \partial_1 h(t_{n}, V_{n})+ \partial_2h(t_{n}, V_{n}) h(t_{n}, V_{n}) \Big)+ \mathcal{O}(\tau^3).
    \end{aligned}
\end{equation}
By subtracting \eqref{C4.11} from \eqref{C4.39}, we obtain
\begin{align}\label{C4.40}
    e_{n+1}= e^{-\tau A} e_{n}+ \tau R_n(V,V_n)+ \delta_{n+1}+\mathcal{O}(\tau^3),
\end{align}
where $R_n$ is given by
\begin{align}
    e^{-\frac{\tau}{2}A} \big( h(t_{n},V)- h(t_{n}, V_{n}) \big)
\end{align}
and similar terms, involving $\partial_1 h$ and $\partial_2 h$. By solving \eqref{C4.40} and using $\| e_0\|=0$, we get
\begin{equation*}
    \Vert e_{n} \Vert \le  \tau \sum_{k=0}^{n-1}  \left\| e^{-\left(n-k-1\right)\tau A}  R_k(V, V_k)\right\| + \sum_{k=0}^{n-1}  \left\| e^{-(n-k-1)\tau A} \left(\delta_{k+1} + C \tau^3\right) \right\|.
\end{equation*}
Using Lemma \ref{Le2}, the local Lipschitz continuity of $f$ and the parabolic smoothing property \eqref{C4.2}, we obtain
\begin{equation}
    \begin{aligned}
    \| e_{n}\|  & \le C \tau \sum_{k=0}^{n-1} \|  e_{k} \| + \sum_{k=0}^{n-2} \| e^{-(n-k-1)\tau A}  A \|  \| \hat{\delta}_{k+1} \| + C\| \delta_{n}\|  + C \tau^2  \\
    &\le   C\tau \sum_{k=0}^{n-1} \Vert e_{k}\Vert + C\tau^3 \sum_{k=1}^{n-1} \frac{1}{k\tau} + C\tau^2\\
    & \le  C\tau \sum_{k=0}^{n-1} \Vert e_{k}\Vert + C\tau^2(1+ \left|  \log \tau \right| ).
    \end{aligned}
\end{equation}
Finally, applying a discrete Gronwall lemma, we get the desired bound \eqref{4.25}. This completes the proof of Theorem \ref{Th: C1}.
\end{proof}

\section{Numerical results in one space dimension}
\label{sec:num-res-1d}
In this section, we present numerical results for the diffusion-reaction problem \eqref{C2.1}, where the diffusion term is represented by the Laplacian, i.e., $D = \partial_{xx}u(t,x)$, and the reaction term is given by $f(t, u(t,x))= u(t,x)^2$ on the domain $\Omega = (0,1)$. For convenience, we refer to the left and right boundary conditions as $b_1$ and $b_2$, respectively. The Laplacian is discretized using a standard second-order centered finite difference scheme with $500$ grid points. To obtain a reference solution, we employ {\sc Matlab} ODE45 with absolute and relative tolerances set to $10^{-9}$. In our simulations, the initial-boundary corrected Strang splitting introduced in Section~\ref{sec:splitting} is referred to as the \emph{IBC Strang splitting}. We evaluate the performance of both classic Strang splitting and IBC Strang splitting when applied to \eqref{C2.1} under various boundary conditions, including Dirichlet, Neumann, Robin, and mixed boundary conditions.

\begin{figure}[ht]
	\centering
	\begin{subfigure}{0.45\textwidth}
		\centering
		\includegraphics[width=\textwidth]{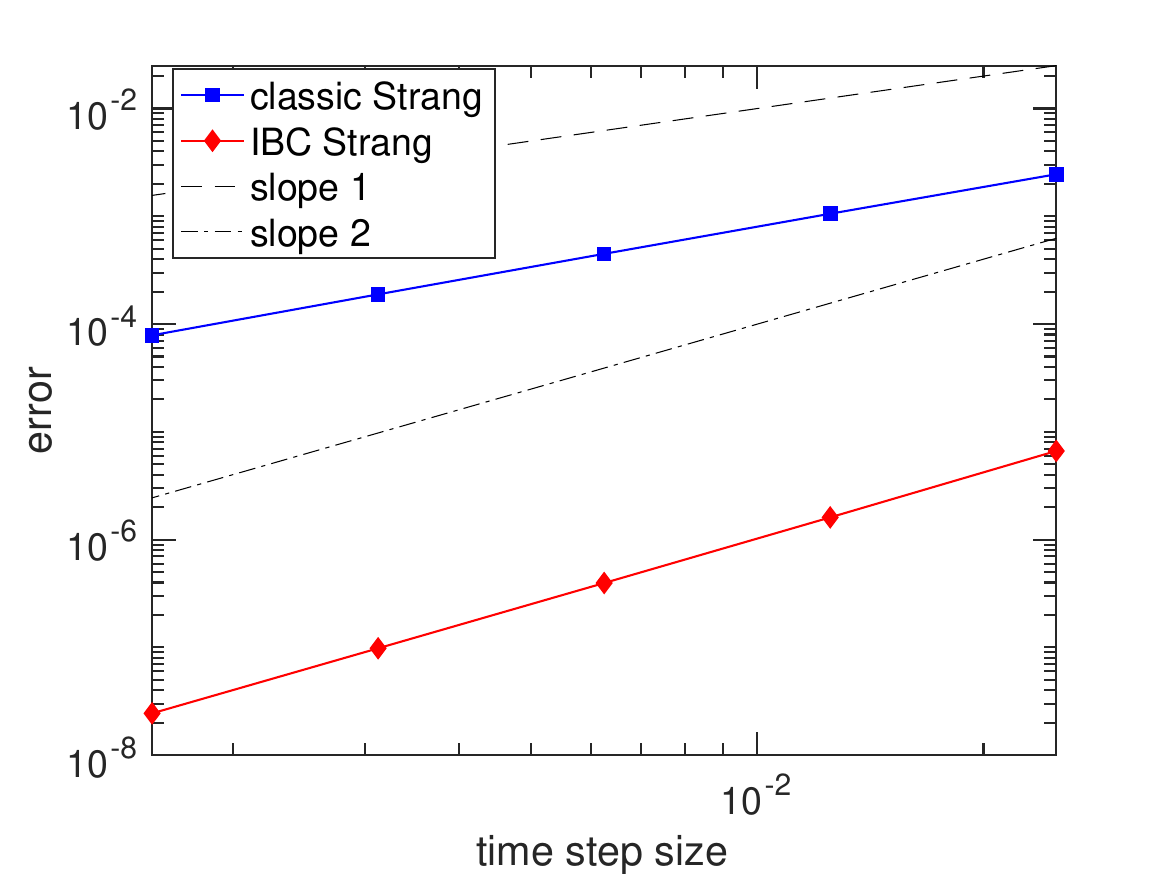}
		\caption{Dirichlet BCs}
		\label{fig:sub1}
	\end{subfigure}
	\hfill
	\begin{subfigure}{0.45\textwidth}
		\centering
		\includegraphics[width=\textwidth]{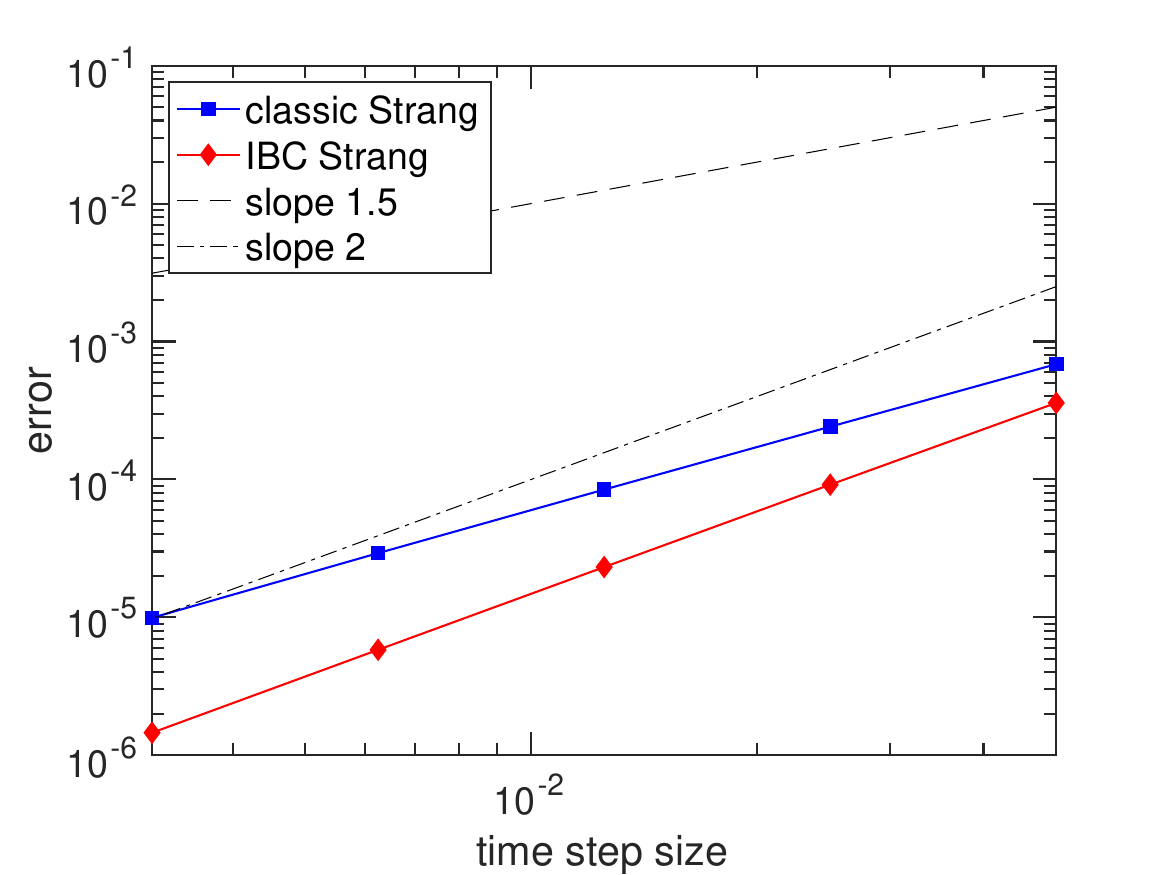}
		\caption{Neumann BCs}
		\label{fig:sub2}
	\end{subfigure}
	
	\vspace{0.5cm}
	
	\begin{subfigure}{0.45\textwidth}
		\centering
		\includegraphics[width=\textwidth]{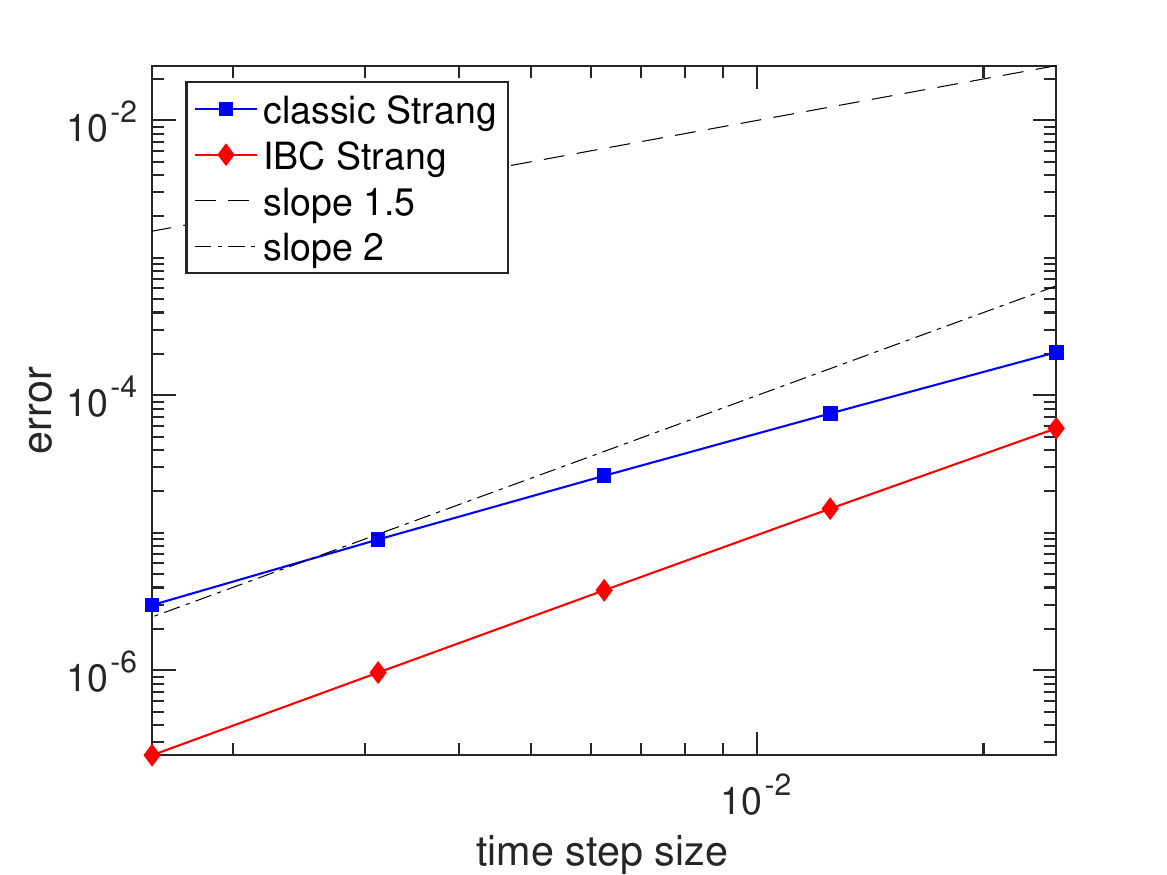}
		\caption{Robin BCs}
		\label{fig:sub3}
	\end{subfigure}
	\hfill
	\begin{subfigure}{0.45\textwidth}
		\centering
		\includegraphics[width=\textwidth]{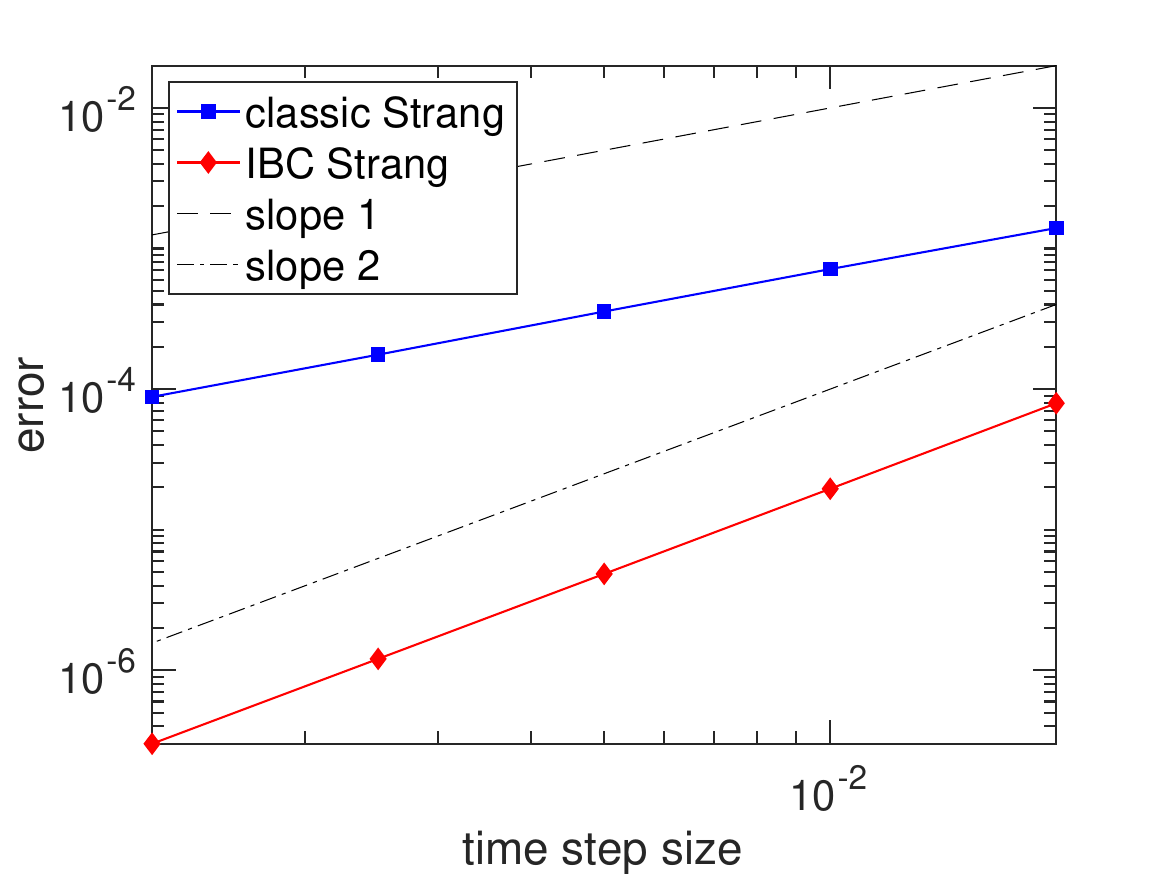}
		\caption{Mixed BCs}
		\label{fig:sub4}
	\end{subfigure}
	
\caption{We solve a one-dimensional diffusion-reaction equation with different types of boundary conditions. The absolute error in the discrete infinity norm is computed at $t=0.5$ by comparing the numerical solution with a reference solution.}
\label{fig:main}
\end{figure}

\begin{exmp}\label{Ex: 5.1}
(Dirichlet boundary conditions, i.e., $\alpha =1$, $ \beta_1 = 0$) In this example, we consider the problem \eqref{C2.1} with the left boundary set to $b_1 = 2$ and the right boundary set to $b_2 = 3$. The initial condition is defined as $u_0 = 2 + \sin\left(\frac{\pi}{2}x \right) $. The corresponding numerical results are presented in Figure \ref{fig:sub1}.
\end{exmp}

\begin{exmp}
(Neumann boundary conditions, i.e., $\alpha =0$, $\beta_1 = 1$) This example is prescribed with $b_1= 1$ and $b_2= 2$. The initial condition is specified as $u_0 = 1 + \frac{2}{\pi} \cos\left( \frac{1}{2} \pi x\right)$ which satisfies the prescribed Neumann boundary conditions. The numerical results are shown in Figure \ref{fig:sub2}.
\end{exmp}

\begin{exmp}
(Robin boundary conditions, i.e., $\alpha = \beta_1 =1$) We consider the problem \eqref{C2.1} with Robin boundary conditions, setting $b_1= 0$ and $b_2 = 1 + \frac{1}{2\pi}$. The initial condition is selected as $u_0(x)= \frac12+ \frac{1}{2\pi}- \frac{1}{2\pi} \cos\left(\pi x \right) $.   Figure \ref{fig:sub3} shows the numerical results.
\end{exmp}

\begin{exmp}\label{Ex: 5.4}
(Mixed boundary conditions)  This example considers the problem \eqref{C2.1} with a Neumann boundary condition $b_1 = 0$ and a Dirichlet boundary condition $b_2 = 2$. We have chosen $u_0(x)= 2- 2  \cos \left( \frac{1}{2}\pi x\right) $ as the initial condition. In Figure \ref{fig:sub4}, we display the numerical results.
\end{exmp}


We observe the expected order reduction to approximately order one for the classic Strang scheme with Dirichlet and mixed boundary conditions (see Figures \ref{fig:sub1} and \ref{fig:sub4}) and reduction to order $1.5$ for Neumann and Robin boundary conditions in the infinity norm (see Figures \ref{fig:sub2} and \ref{fig:sub3}). In contrast, the IBC Strang splitting method shows second-order convergence regardless of the boundary conditions chosen and achieves significantly improved accuracy even for large time step sizes.

\section{Numerical results in two space dimensions}
\label{sec:num-res-2d}
\begin{figure}[ht]
	\centering
	\begin{subfigure}{0.45\textwidth}
		\centering
		\includegraphics[width=\textwidth]{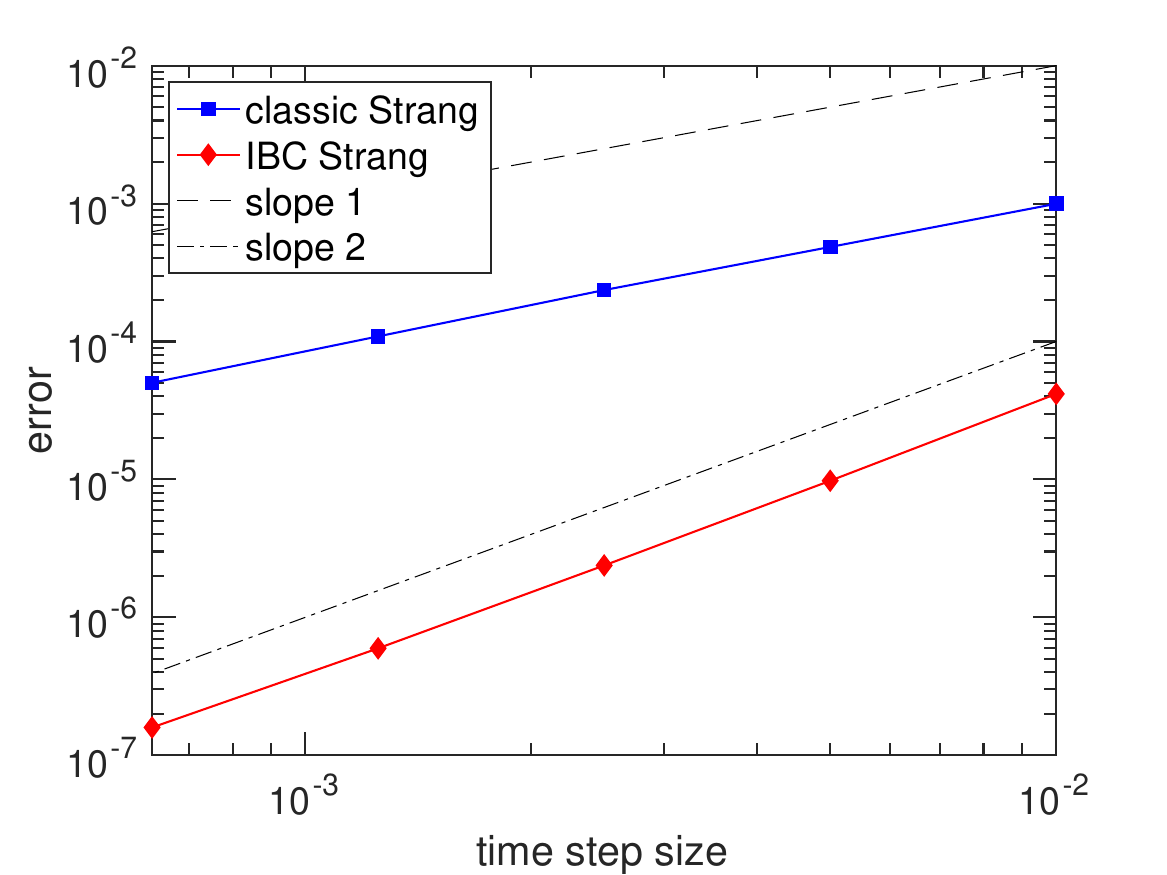}
		\caption{Dirichlet BCs}
		\label{fig:2a}
	\end{subfigure}
	\hfill
	\begin{subfigure}{0.45\textwidth}
		\centering
		\includegraphics[width=\textwidth]{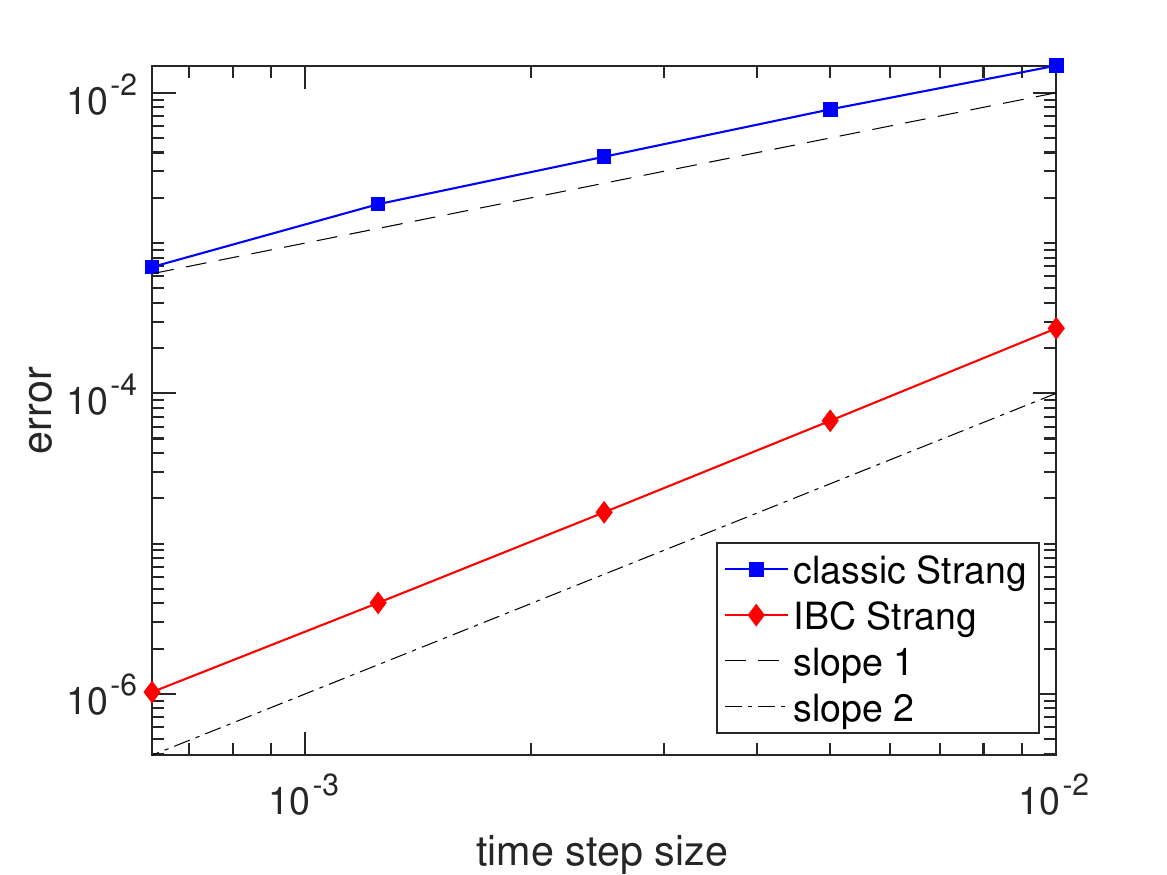}
		\caption{Mixed BCs}
		\label{fig:2b}
	\end{subfigure}
		\caption{We solve a two-dimensional diffusion-reaction equation with different types of boundary conditions. The absolute error in the discrete infinity norm is computed at $t=0.1$ by comparing the numerical solution to a reference solution.}
	\label{fig2}
\end{figure}

In this section, we consider the diffusion-reaction problem \eqref{C2.1} with the nonlinear reaction term given by  $f(t, u(t,x,y)) = u(t,x,y)^2$ on $\Omega = (0,1)^2$, where $D$ is the standard second-order finite difference approximation of the Laplacian, using $50 \times 50$ grid points. The error in the discrete infinity norm is assessed at $t=0.1$ by comparing the numerical solution to a reference solution. This reference solution is obtained using {\sc Matlab} ODE45, with both the absolute and relative tolerances set to $10^{-9}$. We illustrate the performance of both the classic Strang splitting and the IBC Strang splitting methods when applied to \eqref{C2.1}, considering Dirichlet and mixed boundary conditions in Examples \ref{Ex: 6.1} and \ref{Ex: 6.2}, respectively. The numerical results are presented in Figures \ref{fig:2a} and \ref{fig:2b}.

\begin{exmp}\label{Ex: 6.1}
In this example, we prescribe the problem \eqref{C2.1} with $u|_{\partial \Omega} =1$. The initial data is selected as $u_0 = 1+ \sin(\pi x) \sin(\pi y)$ to ensure compliance with the boundary conditions.	
\end{exmp}

\begin{exmp} \label{Ex: 6.2}
For this example, we investigate the problem described by \eqref{C2.1}, applying Dirichlet boundary conditions at the top and bottom and Neumann boundary conditions at the left and right edges. The boundary data $b$
are selected to ensure that the initial condition $u_0(x,y)$ aligns with the specified mixed boundary conditions. The initial value is specified as follows:
	\begin{align*}
		u_0(x,y) = 3+ e^{-10\left( y-\frac{1}{2}\right)^2 } \cos(2 \pi (x+y)).
	\end{align*}
\end{exmp}

In line with the one-dimensional case (see Example \ref{Ex: 5.1} and Example \ref{Ex: 5.4}), it confirms that the classic Strang splitting reduces to a first-order method, while the IBC Strang splitting clearly shows second-order convergence.

\section{Conclusion}
We proposed an initial-boundary corrected splitting method that effectively avoids the order reduction commonly observed when splitting methods are applied to diffusion-reaction problems with non-trivial boundary conditions. In several practical applications, such as those arising from combustion problems, Neumann boundary conditions are of particular interest. Although, modified splitting schemes that restore second order convergence are available in the literature, these often require relatively invasive modifications \cite{doi:10.1137/16M1056250}. In contrast, our new splitting approach can be broadly applied to different boundary conditions without the need to recalculate the modifications, thus allowing for wider applicability.

\section*{Acknowledgements}
This project has received funding from the European Union’s Horizon 2020 research and innovation programme under the Marie Sk\l{}odowska-Curie grant agreement No 847476.


\begin{thebibliography}{10}

\bibitem{10.1016/j.cam.2019.02.023}
I.~Alonso-Mallo, B.~Cano, and N.~Reguera.
\newblock Avoiding order reduction when integrating reaction--diffusion
  boundary value problems with exponential splitting methods.
\newblock {\em Journal of Computational and Applied Mathematics}, 357:228--250,
  2019.

\bibitem{doi:10.1137/19M1257081}
G.~Bertoli and G.~Vilmart.
\newblock Strang splitting method for semilinear parabolic problems with
  inhomogeneous boundary conditions: A correction based on the flow of the
  nonlinearity.
\newblock {\em SIAM Journal on Scientific Computing}, 42(3):A1913--A1934, 2020.

\bibitem{f632bc39-61b5-3c7b-90cc-3a08f6b7fde1}
S.~Descombes.
\newblock Convergence of a splitting method of high order for
  reaction-diffusion systems.
\newblock {\em Mathematics of Computation}, 70(236):1481--1501, 2001.

\bibitem{EINKEMMER201876}
L.~Einkemmer, M.~Moccaldi, and A.~Ostermann.
\newblock Efficient boundary corrected {S}trang splitting.
\newblock {\em Applied Mathematics and Computation}, 332:76--89, 2018.

\bibitem{doi:10.1137/140994204}
L.~Einkemmer and A.~Ostermann.
\newblock Overcoming order reduction in diffusion-reaction splitting. {P}art 1:
  Dirichlet boundary conditions.
\newblock {\em SIAM Journal on Scientific Computing}, 37(3):A1577--A1592, 2015.

\bibitem{doi:10.1137/16M1056250}
L.~Einkemmer and A.~Ostermann.
\newblock Overcoming order reduction in diffusion-reaction splitting. {P}art 2:
  Oblique boundary conditions.
\newblock {\em SIAM Journal on Scientific Computing}, 38(6):A3741--A3757, 2016.

\bibitem{Einkemmer2016ACO}
L.~Einkemmer and A.~Ostermann.
\newblock A comparison of boundary correction methods for {S}trang splitting.
\newblock {\em Discrete and Continuous Dynamical Systems - B},
  23(7):2641--2660, 2016.

\bibitem{GERISCH2002159}
A.~Gerisch and J.~G. Verwer.
\newblock Operator splitting and approximate factorization for
  taxis--diffusion--reaction models.
\newblock {\em Applied Numerical Mathematics}, 42(1):159--176, 2002.

\bibitem{hairer2013geometric}
E.~Hairer, C.~Lubich, and G.~Wanner.
\newblock {\em Geometric Numerical Integration: Structure-Preserving Algorithms
  for Ordinary Differential Equations}.
\newblock Springer Series in Computational Mathematics. Springer, Berlin
  Heidelberg, 2013.

\bibitem{HANSEN20121428}
E.~Hansen, F.~Kramer, and A.~Ostermann.
\newblock A second-order positivity preserving scheme for semilinear parabolic
  problems.
\newblock {\em Applied Numerical Mathematics}, 62(10):1428--1435, 2012.

\bibitem{henry1981geometric}
D.~Henry.
\newblock {\em Geometric Theory of Semilinear Parabolic Equations}.
\newblock Lecture notes in mathematics. Springer, Berlin Heidelberg New York,
  1981.

\bibitem{HUNDSDORFER1995191}
W.~Hundsdorfer and J.~G. Verwer.
\newblock A note on splitting errors for advection-reaction equations.
\newblock {\em Applied Numerical Mathematics}, 18(1):191--199, 1995.

\bibitem{hundsdorfer2013numerical}
W.~Hundsdorfer and J.~G. Verwer.
\newblock {\em Numerical Solution of Time-Dependent
  Advection-Diffusion-Reaction Equations}.
\newblock Springer Series in Computational Mathematics. Springer, Berlin
  Heidelberg, 2013.

\bibitem{MCKENNEY1995348}
A.~McKenney, L.~Greengard, and A.~Mayo.
\newblock A fast {P}oisson solver for complex geometries.
\newblock {\em Journal of Computational Physics}, 118(2):348--355, 1995.

\bibitem{mclachlan_quispel_2002}
R.~I. McLachlan and G.~R.~W. Quispel.
\newblock Splitting methods.
\newblock {\em Acta Numerica}, 11:341--434, 2002.

\bibitem{pazy2012semigroups}
A.~Pazy.
\newblock {\em Semigroups of Linear Operators and Applications to Partial
  Differential Equations}.
\newblock Applied Mathematical Sciences. Springer, New York, 2012.

\bibitem{doi:10.1137/0705041}
G.~Strang.
\newblock On the construction and comparison of difference schemes.
\newblock {\em SIAM Journal on Numerical Analysis}, 5(3):506--517, 1968.

\end{thebibliography}

\end{document}